\documentclass[12pt]{amsart}
\usepackage{amssymb}
\usepackage{a4wide}
\usepackage{color}
 \usepackage{boxedminipage}
\usepackage{enumerate, graphicx}
\usepackage{todonotes}
\usepackage{subfig}
\usepackage{wrapfig}
\usepackage{cite}
\usepackage{caption}
\usepackage[colorlinks,citecolor=blue,linkcolor=blue]{hyperref}
\usepackage[capitalise,nameinlink]{cleveref}

  \makeatletter
  \CheckCommand*\refstepcounter[1]{\stepcounter{#1}%
      \protected@edef\@currentlabel
       {\csname p@#1\endcsname\csname the#1\endcsname}%
  }
  \renewcommand*\refstepcounter[1]{\stepcounter{#1}%
    \protected@edef\@currentlabel
      {\csname p@#1\expandafter\endcsname\csname the#1\endcsname}%
  }
  \def\labelformat#1{\expandafter\def\csname p@#1\endcsname##1}
  \DeclareRobustCommand\Ref[1]{\protected@edef\@tempa{\ref{#1}}%
     \expandafter\MakeUppercase\@tempa
  }
  \makeatother

  \makeatletter
  \newcommand{\numberlike}[2]{%
     \expandafter\def\csname c@#1\endcsname{%
         \expandafter\csname c@#2\endcsname}%
  }
  \makeatother

  \def\DefaultNumberTheoremWithin{section}

  \theoremstyle{plain}
  \newtheorem{Lemma}{Lemma}
     \numberwithin{Lemma}{\DefaultNumberTheoremWithin}
     \labelformat{Lemma}{Lemma~#1}
  
     \numberwithin{Claim}{\DefaultNumberTheoremWithin}
     \numberlike{Claim}{Lemma}
     \labelformat{Claim}{Claim~#1}

  \newtheorem{Theorem}{Theorem}
     \numberwithin{Theorem}{\DefaultNumberTheoremWithin}
     \numberlike{Theorem}{Lemma}
     \labelformat{Theorem}{Theorem~#1}
  \newtheorem{Corollary}{Corollary}
     \numberwithin{Corollary}{\DefaultNumberTheoremWithin}
     \numberlike{Corollary}{Lemma}
     \labelformat{Corollary}{Corollary~#1}
  
     \numberwithin{Proposition}{\DefaultNumberTheoremWithin}
     \numberlike{Proposition}{Lemma}
     \labelformat{Proposition}{Proposition~#1}
  
     \numberwithin{Conjecture}{\DefaultNumberTheoremWithin}
     \numberlike{Conjecture}{Lemma}
     \labelformat{Conjecture}{Conjecture~#1}
  
     \numberwithin{Situation}{\DefaultNumberTheoremWithin}
     \numberlike{Situation}{Lemma}
     \labelformat{Situation}{Situation~#1}
 
     \numberwithin{Note}{\DefaultNumberTheoremWithin}
     \numberlike{Note}{Lemma}
     \labelformat{Note}{Note~#1}
     
  \theoremstyle{definition}
  \newtheorem{Definition}{Definition}
     \numberwithin{Definition}{\DefaultNumberTheoremWithin}
     \numberlike{Definition}{Lemma}
     \labelformat{Definition}{Definition~#1}

  \theoremstyle{definition}
  
     \numberwithin{Question}{\DefaultNumberTheoremWithin}
     \numberlike{Question}{Lemma}
     \labelformat{Question}{Question~#1}

  \theoremstyle{definition}
  
     \numberwithin{Problem}{\DefaultNumberTheoremWithin}
     \numberlike{Problem}{Lemma}
     \labelformat{Problem}{Problem~#1}

     \theoremstyle{remark} 
     \numberwithin{Remark}{\DefaultNumberTheoremWithin}
     \numberlike{Remark}{Lemma}
     \labelformat{Remark}{Remark~#1}
  \theoremstyle{remark}

  \newtheorem{Example}{Example}
     \numberwithin{Example}{\DefaultNumberTheoremWithin}
     \numberlike{Example}{Lemma}
     \labelformat{Example}{Example~#1}
  
     \labelformat{Case}{Case~#1}
     \numberwithin{Case}{Lemma}
  
     \labelformat{Step}{Step~#1}
     \numberwithin{Step}{Lemma}

  \newtheorem*{Example*}{Example}

\labelformat{section}{Section~#1}
  \labelformat{subsection}{Section~#1}
  \labelformat{subsubsection}{Section~#1}

\title{The Chain Matrix of Bouquets of Geometric Lattices and its Determinant}
\author{Winfried Hochst\"attler}
    \address{FernUniversit\"at in Hagen \\ 
          Fakult\"at f\"ur Mathematik und Informatik \\
          58084 Hagen\\
          Germany}
     \email{winfried.hochstaettler@fernuni-hagen.de}
     
\author{Sophia Keip}
    \address{FernUniversit\"at in Hagen \\ 
          Fakult\"at f\"ur Mathematik und Informatik \\
          58084 Hagen\\
          Germany}
     \email{sophia.keip@fernuni-hagen.de}

\begin{document}
\begin{abstract}
This work builds on Varchenko et al's introduction of bilinear forms for hyperplane arrangements, where the determinant of the associated matrices factorizes into simple components. While one of the determinant formula developed by Varchenko has been generalized to complexes of oriented matroids (COMs) already, this question was open for another, distinct form. Motivated by work from Varchenko and Brylawski, who generalized the alternative bilinear form and its determinant formula from hyperplane arrangements to matroids, we examine whether this formula can similarly be generalized to COMs. Our findings affirm this generalization, and we further extend the determinant formula to bouquets of geometric lattices as introduced by Laurent et al.
\end{abstract}

\maketitle
\section{Introduction}
\subsection{Background and Motivation}
Varchenko and coauthors introduced various bili\-near forms for hyperplane arrangements. A common feature of these forms is that the determinant of the matrix representing them can be factorized into a product with remarkably simple factors. The determinant formula in \cite{varchenko1995multidimensional} has been generalized to oriented matroids (OMs) \cite{hochstattler2019varchenko}. Building on this, we extended the formula further to complexes of oriented matroids (COMs), which generalize OMs, in \cite{hochstattler2022signed}. During this work, Varchenko pointed us to his related work with Brylawski \cite{brylawski1997determinant}, which explores a distinct bilinear form with a determinant formula for matroids. The bilinear form in Varchenko and Brylawski's work was initially introduced by Schechtman and Varchenko for complex hyperplane arrangements \cite{schechtman1991arrangements}. While the matrix in \cite{varchenko1995multidimensional} is indexed by the full-dimensional cells of a hyperplane arrangement, in \cite{brylawski1997determinant} it is defined on the flag space of the arrangement. For a COM, this space corresponds to the zero sets of its covectors. Since the determinant formula in \cite{varchenko1995multidimensional} could be extended to COMs, a natural question arose as to whether the same generalization would apply to \cite{brylawski1997determinant}. In this work, we affirmatively answer that question. Furthermore, we achieve a more general result by proving the determinant formula for bouquets of geometric lattices introduced by Laurent et al. in \cite{laurent1989bouquets}. In our generalization, we work with maximal chains, which may have different lengths instead of the classical flags of geometric lattices. Therefore, we refer to the matrix representing the bilinear form as the chain matrix.

\subsection{Bouquets of Geometric Lattices}
We begin with some basics about lattices, following \cite{wachs2006poset}. Let $P$ be a partially ordered set (a \emph{poset}). For $x, y \in P$, the \emph{meet} of $x$ and $y$, denoted $x \land y$, is a unique element in $P$ that is less than or equal to both $x$ and $y$ and is greater than any other element with this property. Similarly, the \emph{join} of $x$ and $y$, denoted $x \lor y$, is a unique element in $P$ that is greater than or equal to both $x$ and $y$ and is smaller than any other element with this property. In simpler terms, having a meet means there exists a greatest lower bound for $x$ and $y$, while having a join means there exists a smallest upper bound for $x$ and $y$. A poset $P$ is called a \emph{meet semilattice} (or \emph{join semilattice}) if every pair of elements in $P$ has a meet (or join). For $x \neq y \in P$ with $x \leq y$, we say that $y$ \emph{covers} $x$ if $x \leq z \leq y$ implies $z \in \{x, y\}$. This is denoted by $x \lessdot y$.

If a poset has a unique minimal element, it is denoted by $\hat{0}$, and if it has a unique maximal element, it is denoted by $\hat{1}$. An element $a \in P$ such that $\hat{0} \lessdot a$ is called an \emph{atom}. Finally, a \emph{chain} $C$ in a poset is a subset of $P$ in which every pair of elements is comparable (i.e., $C$ is totally ordered).

To understand bouquets of geometric lattices, we need to define geometric lattices first. Geometric lattices are a special class of lattices that arise from the lattice of flats of a matroid, making them closely related to combinatorial geometry.

\begin{Definition}[(geometric) Lattice]
A poset $P$ is called a \emph{lattice}, if every two elements $x,y \in P$ have a meet and a join. A lattice is called \emph{geometric} if it is
\begin{itemize}
\item[1.]\emph{atomic}, i.e.
\begin{align*}
\forall x \in P:\,x = \bigvee\limits_{\hat{0}\lessdot a\leq x} a
\end{align*} 
\item[2.]\emph{semimodular}, i.e.
\begin{align*}
\forall x,y \in P: \left((x\land y) \lessdot x \Rightarrow y \lessdot(x\lor y)\right).
\end{align*} 
\end{itemize}
\end{Definition}

\begin{figure}
{\includegraphics[scale=1]{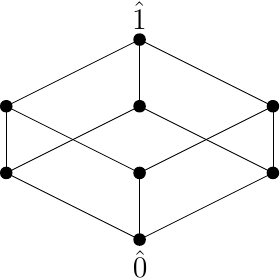}}\hspace{2cm}
{\includegraphics[scale=1]{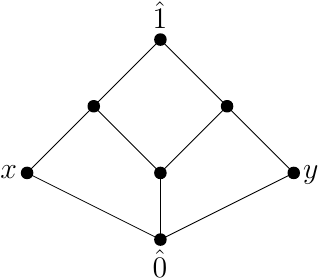}}
\caption{Examples of a lattices with a $\hat{1}$ and a $\hat{0}$ element. The lattice to the left is geometric. The lattice to the right is not geometric, since semimodularity is not fulfilled for $x$ and $y$. It is neither a bouquet of geometric lattices.\label{fig:geom}}
\end{figure}

This definition is illustrated in Figure \ref{fig:geom}. Bouquets of geometric lattices extend the classical concept of geometric lattices:

\begin{Definition}[\cite{laurent1989bouquets} Bouquets of geometric Lattices]\label{def:bouquetlattice}
A poset $P$ is a bouquet of geometric lattices if $P$ is a meet semilattice in which every interval is a geometric lattice.
\end{Definition}

\begin{figure}
\includegraphics[scale=1]{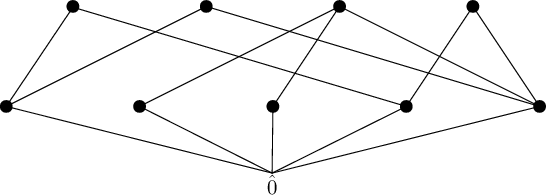}
\caption{Example of a bouquet of geometric lattices.\label{fig:bouquetgeom}}
\end{figure}

An example is shown in Figure \ref{fig:bouquetgeom}. Just as geometric lattices are linked to matroids, bouquets of geometric lattices are closely associated with simple bouquets of matroids. These connections will be explored further in \ref{sec:bouquetsMatr}. Furthermore, bouquets of geometric lattices generalize the framework of geometric semilattices that have been studied by Wachs and Walker \cite{wachs1985geometric}.

 A bouquet of geometric lattices is ranked with a rank function $r$. This means every unrefinable chain from a minimal element to a fixed element $x$ has the same length $r(x)$. The \emph{M\"obius function} $\mu$ of a poset is defined by 
\begin{align*}
& \mu(x,y)=0   \text{ for all }x \not \le y \in \mathcal{P}\\
&\mu(x,x) = 1\text{ for all }x \in \mathcal{P}\\
&\mu(x,y) = -\sum_{x \leq z <y} \mu(x,z)\text{ for all }x <y \in \mathcal{P}.  
\end{align*}
For $\mu(\hat0,x)$ we write $\mu(x)$ and we denote the unsigned version by $\mu^+(x,y) = |\mu(x,y)|$. Crapo's \emph{beta function} $\beta$ \cite{crapo1967higher} of a poset element is defined by
\begin{align*}
    \beta(x) = (-1)^{r(x)}\sum_{y:y\leq x} \mu(y) r(y).
\end{align*}

\medskip

We now have all the necessary components to define the \emph{cumulated rho function}, a key element in our factorization formula. Let $r_1, \dots, r_k$ denote the maximal elements of $P$ such that $x \leq r_i$ for $i = 1, \dots, k$. The cumulated rho function $\rho$ is then defined by 
\begin{align*}
    \rho_P(x) = \beta(x)\sum_{i = 1}^k\mu^+(x,r_i).
\end{align*}

\medskip

In the introduction of \cite{brylawski2000mobius}, Brylawski gives several interpretations of this function.
We will see a application of those formulas later in \ref{ex:formula}, when we demonstrate our factorization formula on a concrete example.

\pagebreak[3]

\subsection{The Chain Matrix}
In order to generalize the matroid bilinear form of Brylawski and Varchenko \cite{brylawski1997determinant}, we need to look at the chains of a bouquet of geometric lattices $P$. If $\{x_1,x_2,\dots,x_k\}$ is a maximal chain of $P$ this means, there is no $y \in P$ such that $x_k \lessdot y$ and
\begin{align*}
C=[x_1\lessdot x_2\lessdot\dots\lessdot x_k].
\end{align*}
We denote $|C|=k$ for the size of $C$.  In the following we only consider maximal chains, therefore we omit the term maximal in following. Let $A_k = (a_1,\dots,a_k)$ be a $k-$tuple of atoms. We say a chain $C=\{x_1,x_2,\dots,x_k\}$ is \emph{generated} by $A_k$, written $ch(A_k) = C$, if 
\begin{align} \label{adj}
x_i = a_1\lor \dots \lor a_i
\end{align}
for all $i = 1,\dots,k$. 
Let us now label every element $x\in P$ with an atom $a \leq x$. We denote this labeling by $l(x)$. We call a chain $C = \{x_1,x_2,\dots,x_k\}$ \emph{neat} , if $l(x_i)\leq x_i$ but $l(x_i) \nleq x_{i-1}$ for all $i$. We call a labeling \emph{convex} if for every $a < x' < x$, such that $l(x) = a$, we have $l(x')=a$. In this case a chain is neat if and only if all of its elements have distinct labels. A very natural convex labeling is the \emph{min-labeling}. Here, we assume a linear order on the atoms $a_1 \prec a_2 \prec \dots \prec a_n$ and define $l(x) = a_i$, where $a_i$ is the smallest atom such that $a_i \leq x$.

\begin{Example}[Min-Labeling]\label{ex:min-labeling}
We look at the example from Figure \ref{fig:bouquetgeom}. Let us assign names to the elements of the poset first:
\begin{center}
  \includegraphics[scale = 1]{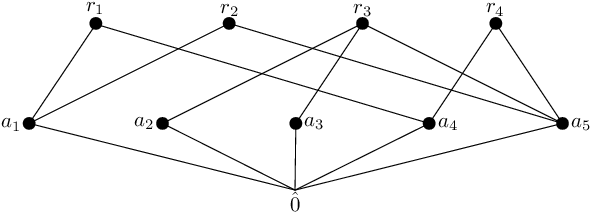}  
  \captionof{figure}{Bouquet of geometric lattices with atoms $a_1, \dots, a_5$.}\label{fig:bouquetgeom1}
\end{center}

With the min-labeling we get $l(a_i) = a_i$, $i=1,\dots,5$ for all atoms and $l(r_1) = a_1$, $l(r_2) = a_1$, $l(r_3) = a_2$ and $l(r_4) = a_4$ for the other elements. The neat chains are $[a_5\lessdot r_4]$, $[a_5\lessdot r_3]$, $[a_5\lessdot r_2]$, $[a_4\lessdot r_1]$ and $[a_3\lessdot r_3]$ since there all elements have distinct labels.
\end{Example}

Consistent with \ref{ex:min-labeling}, let $r_1,\dots,r_n$ be the maximal elements of $P$. We call the set of neat maximal chains that have $r_i$ as biggest element the \emph{neat chain family} $C(r_i)$ of $r_i$. Furthermore, let $\mathcal{C}= \bigcup_{i=1}^n C(r_i)$ be the set of all neat chains of $P$.

\begin{Example*}[1.3 continued]
For the bouquet of geometric lattices of \ref{ex:min-labeling} the maximal elements are $r_1,\dots,r_4$. We therefore have $C(r_1) = \{[a_4\lessdot r_1]\}$, $C(r_2) = \{[a_5\lessdot r_2]\}$, $C(r_3)= \{[a_3\lessdot r_3],[a_5\lessdot r_3]\}$ and $C(r_4) = \{[a_5\lessdot r_4]\}$.
\end{Example*}

We now assign variables from the ring $\mathbb{Z}[w_1,\dots,w_n]$ to the atoms of our bouquet of geometric lattices, considered as weights. The \emph{weight} of an element $x \in P$ is defined by
\begin{align*}
w(x) = \sum_{a_i \leq x}w_i.
\end{align*}
We are now ready to define the chain matrix.

\begin{Definition}[Chain Matrix]\label{def:BVmatrix}
Let $P$ be a bouquet of geometric lattices and $\mathcal{C}$ its neat chains. The chain matrix $\mathfrak{C}$ is a symmetric $\mathcal{C}\times \mathcal{C}-$Matrix defined by
\begin{align*}
\mathfrak{C}_{C,C'} = \sum\limits_{ch(a_{i_1},\dots,a_{i_k}) = C,\, ch(a_{\sigma(i_1)},\dots,a_{\sigma(i_k)}) = C'} \text{sgn} (\sigma)\cdot w_{i_1} \dots  w_{i_k} 
\end{align*}
where $C,C' \in \mathcal{C}$ and $\sigma$ is a permutation of $\{i_1,\dots,i_k\}$, $k= |C|= |C'|$. 
\end{Definition}

\begin{Example*}[1.3 continued]\label{ex:chainmatrixgeomlattice}
The chain matrix of the bouquet of geometric lattices from \ref{ex:min-labeling} is
\begin{align*}
\bordermatrix{
  & [a_4\lessdot r_1]   & [a_5\lessdot r_2]   & [a_5\lessdot r_3]   & [a_3\lessdot r_3] & [a_5\lessdot r_4]  \cr
[a_4\lessdot r_1] & w_1\cdot w_4 & 0 & 0 & 0 & 0 \cr
[a_5\lessdot r_2] & 0 & w_1\cdot w_5 & 0 & 0 & 0\cr
[a_5\lessdot r_3] & 0 & 0 & w_2 \cdot w_5 + w_3\cdot w_5  & -w_3\cdot w_5 & 0\cr
[a_3\lessdot r_3] & 0 & 0 & -w_3\cdot w_5  & w_2 \cdot w_3 + w_3 \cdot w_5 & 0\cr
[a_5\lessdot r_4] & 0 & 0 & 0 & 0 & w_4 \cdot w_5\cr
}
\end{align*}

\vspace{2mm}
We discuss the entries of the third row for a better understanding. We start with the third entry where $C = C' = [a_5\lessdot r_3]$. $(a_5, a_2)$ and $(a_5,a_3)$ are generators for $[a_5\lessdot r_3]$, so we have the two summands $w_2 \cdot w_5$ and $w_3\cdot w_5$. Since $C = C'$, $\text{sgn}(\sigma) = 1$, so both summands have a positive sign. For the fourth entry we see that $(a_3,a_5)$ generates $[a_3\lessdot r_3]$ and $(a_5,a_3)$ generates $[a_5\lessdot r_3]$. One sees that $\text{sgn}(\sigma) = -1$, therefore we get the entry $-w_3 \cdot w_5$.
The first, second and fifth entry is zero, since the corresponding chains do not share generating tuples of atoms.
\end{Example*}
\subsection{Results}
We will prove the following factorization formula

\begin{Theorem}[Chain Matrix Determinant for Bouquets of geometric Lattices]\label{main}
Let $P$ be a bouquet of geometric lattices and $\mathfrak{C}$ its chain matrix. Then 
\begin{align*}
det(\mathfrak{C}) = \pm \prod\limits_{x \in P} w(x)^{\rho_P(x)}.
\end{align*}
\end{Theorem}
This is a generalization of Theorem 4.16 of \cite{brylawski1997determinant}. Furthermore, we will link complexes of oriented matroids to bouquets of geometric lattices in \ref{sec:COMs}, leading to the generalization of Brylawski and Varchenkos theorem for them as well as corrollary of our theorem.

\subsection{Related Work} Very recently, in \cite{eberhardt2024intersection}, Eberhardt and Mautner introduce an integer-valued intersection matrix for real affine hyperplane arrangements and, more generally, for affine oriented matroids, proving a determinant formula for it. Similar to our work, their bilinear form is closely related to the bilinear form introduced by Brylawski and Varchenko \cite{brylawski1997determinant}, and they also utilize the determinant theorem established there to prove their formula. Exploring the connection between their results and ours would be highly intriguing, particularly in determining whether their framework could also be to COMs.

\section{Proof of \ref{main}}
In this section we will first explain the result of Brylawski and Varchenko \cite{brylawski1997determinant} for matroids. In the second part we will use this result to prove \ref{main}.

\subsection{The Determinant of the Chain Matrix of a Matroid}
One special case of bouquets of geometric lattices are the flat lattices of matroids, which are even geometric lattices, as indicated in the introduction. In order to look at the determinant of the chain matrix of that special case, let us start with the definition.
\begin{Definition}[Matroid]
A finite groundset $E$ together with a collection $\mathcal{I}$ of subsets of $E$ is called matroid $\mathcal{M}=(E,\mathcal{I})$ if
\begin{enumerate}
\item $\emptyset \in \mathcal{I}$,
\item $I_1 \subset I_2 \in \mathcal{I}\;\Rightarrow\; I_1 \in \mathcal{I}$,
\item $I_1, I_2 \in \mathcal{I},\; |I_1|<|I_2|\;\Rightarrow\; \exists\, e\in I_2\backslash I_1:\; I_1\cup e \in \mathcal{I}$.
\end{enumerate}
The elements of $\mathcal{I}$ are called \emph{independent sets}. A subset of $E$ which is not in $\mathcal{I}$ is called dependent. 
\end{Definition}
A \emph{circuit} of a matroid is a minimal dependent subset $S \subseteq E$, i.e. all proper subsets of $S$ are independent. A matroid is called \emph{simple} if it has no circuits of cardinality one or two.
Let $S \subseteq E$. The \emph{rank} $r(S)$ of $S$ is defined by the size of the largest independent set contained in $S$. The rank of the whole matroid $r(\mathcal{M})$ is $r(E)$. $S \subseteq E$ is called \emph{flat} if for all $e\in E\backslash S$ $r(S\cup e)>r(S)$. We define a partial order on the flats of a matroid where $S<S'$ if $S \subset S'$. 


The poset defined on the flats of a matroid forms a geometric lattice. In fact, there is a bijection between simple matroids and geometric lattices. For the flat lattice of a matroid, the chain matrix corresponds to the \emph{bilinear form of a matroid} as introduced in \cite{brylawski1997determinant}. The chains of the flat lattice of a matroid are referred to as \emph{flags}, and in \cite{brylawski2000mobius}, the bilinear form of a matroid is also called the \emph{flag matrix}. Varchenko and Brylawski proved the following factorization formula for the determinant of this matrix.

\begin{Theorem}[Brylawski-Varchenko \cite{brylawski1997determinant}]\label{BV}
Let $\mathcal{M}$ be a matroid and $\mathfrak{F}$ its flag matrix. Then 
\begin{align*}
det(\mathfrak{F}) = \pm \prod\limits_{K \text{ flat of } \mathcal{M}} w(K)^{b_K}
\end{align*}
where $w(K) = \sum_{e \in K}w_e$ and $b_K = \beta(K)|\mu(K,\hat1)|$.
\end{Theorem}

To see that $b_K$ this is the accumulated rho function we defined earlier, note that the flat lattice of a matroid has a unique maximal element. 
Our goal is now to use this determinant formula to prove the formula for bouquets of geometric lattices.

\subsection{The Chain Matrix Determinant of a Bouquet of geometric lattices}
To prove theorem \ref{main} we need the following lemma.

\begin{Lemma}\label{nullentries} Let $P$ be a bouquet of geometric lattices and $\mathcal{C}=C(r_1) \cup C(r_2) \cup \dots \cup C(r_n)$ its neat chains. Let $C \in C(r_i)$ and $C' \in C(r_j)$, $i\neq j$. Then $\mathfrak{C}_{C,C'}=0$.
\end{Lemma}

\begin{proof}
If $|C| \neq |C'|$ there is nothing to show, so let us assume $|C| = |C'| = k$. Let $A$ be the set of atoms of $P$. We need to show that the set 
\begin{align*}
\{(a_{i_1},\dots,a_{i_k})\subseteq A:\; ch(a_{i_1},\dots,a_{i_k}) = C, ch(a_{\sigma(i_1)},\dots,a_{\sigma(i_k)}) = C'\}
\end{align*}
is empty. If $ch(a_{i_1},\dots,a_{i_k}) = C$ and $ch(a_{\sigma(i_1)},\dots,a_{\sigma(i_k)}) = C'$, then because of (\ref{adj}) we have $r_i = a_{i_1}\lor \dots \lor a_{i_k} = r_j$, since taking a join does not depend on the order of the elements. This is a contradiction to the assumption $i \neq j$.
\end{proof}

Now we can prove \ref{main}.

\begin{proof}
Let $P$ be a bouquet of geometric lattices and $\mathfrak{C}$ its chain matrix. Let $\mathcal{C}=C(r_1) \cup C(r_2) \cup \dots \cup C(r_n)$ be the neat chains of $P$. We now reorder the rows and columns of $\mathfrak{C}$ such that columns/rows belonging to the chains of $C(r_i)$ are one after the other for all $i=1,\dots,n$. Since we swap an equal amount of columns and rows here, the sign of the determinant does not change. From \ref{nullentries} we know, that an entry of the matrix is only non zero, if it is indexed by chains belonging to the same $C(r_i)$. Since we ordered the rows and columns properly, we have a block diagonal matrix, where each block $B_{r_i}$ is indexed by the chains belonging to only one $C(r_i)$. Since the interval $[\emptyset,r_i]$ is geometric, it is the flat lattice for some matroid. So we can compute the determinant of $B_{r_i}$ by \cite{brylawski1997determinant}.
\begin{align*}
det(B_{r_i}) = \pm \prod\limits_{x \leq r_i} w(x)^{\beta(x)|\mu(x,\hat1)|}.
\end{align*}

Since the determinant of the whole matrix is the product of the determinants of the blocks on its diagonal, our claim follows.
\end{proof}

\begin{Example*}[1.3 continued] \label{ex:formula}
    We continue with example \ref{ex:chainmatrixgeomlattice}. Before giving the whole factorization let us derive the exponent of some elements step by step. We start with the rank 1 element $a_1$. $\mu(a_1) = -1$ and therefore $\beta(x) = (-1)^1 \cdot (-1) \cdot 1 = 1$. There are two maximal elements that are bigger than $a_1$: $r_1, r_2$. $\mu^+(a_1,r_i) = 1$, $i=1,2$, so in total we get $\rho_P(a_1) = 1 \cdot (1+1) = 2$. Now look at $r_1$. This element has rank 2, $\mu(r_1) = \mu^+(r_1) = 1$ and $\beta(r_1) = (-1)^2(-1+-1+2) = 0$, so $\rho_P(r_1) = 0$. The other exponents can be obtained in the same way and in total we get
    \begin{align*}
        det(\mathfrak{C}) = w_5^3 \cdot w_4^2 \cdot w_3 \cdot w_2 \cdot w_1^2 \cdot (w_2 + w_3 + w_5)
    \end{align*}
    which is in line with \ref{main}.
\end{Example*}

\section{Applications}
\subsection{Zero Sets of Complexes of oriented Matroids}\label{sec:COMs}
One application of our theorem can be found in the context of \emph{complexes of oriented matroids (COMs)}. This combinatorial structure consists of a finite groundset together with vectors containing $+, -$ and $0$ entries, so called sign vectors and was introduced as a generalization of lopsided sets and oriented matroids \cite{bandelt2018coms}. We will show that the poset on the zero sets of COMs induced by inclusion forms a bouquet of geometric lattices. Before we show that, let us introduce COMs formally. A COM can be defined by a finited groundset $E$ together with a collection of sign vectors $\mathcal{L} \subseteq \{0,+,-\}^{|E|}$. Before we can formulate the axioms defining a COM, we need to define the following terms.
\begin{Definition} Let $E$ be a finite set and $X,Y \in \{0,+,-\}^{|E|}$. We call
\begin{align*}
    S(X,Y) = \{e \in E: X_e = -Y_e \neq 0\}. 
\end{align*}
the \emph{separator} of $X$ and $Y$. The \emph{composition} of $X$ and $Y$ is defined by
\begin{align*}
(X \circ Y)_e = \begin{cases}
                        X_e & \text{ if } X_e \neq 0,\\
                        Y_e & \text{ if } X_e = 0\\
                       \end{cases} 
                    \forall e\in E.
\end{align*} 
The \emph{zero set} of $X$ is denoted by
\begin{align*}
z(X) = E\backslash\underline{X} = \{e \in E: X_e =  0\},
\end{align*}
where $\underline{X}$ denotes the \emph{support} of $X$
\begin{align*}
\underline{X} = \{e \in E: X_e \neq  0\}.
\end{align*}
\end{Definition}

Now we can define the term COM.
\begin{Definition}[Complex of Oriented Matroids (COM)] Let $E$ be a finite set and $\mathcal{L} \subseteq \{0,+,-\}^{|E|}$. The pair $\mathcal{M}=(E,\mathcal{L})$ is called a COM, if $\mathcal{L}$ satisfies 
\begin{itemize}
    \item[(FS)]  Face Symmetry 
    \begin{align*}
\forall X,Y \in \mathcal{L}: X \circ (-Y) \in \mathcal{L}.
\end{align*}
\end{itemize}
and
\begin{itemize}
\item[(SE)] Strong Elimination 
\begin{align*}
&\forall X,Y \in \mathcal{L}\, \forall e \in S(X,Y)\, \exists Z \in \mathcal{L}: \\
&Z_e=0 \text{ and }\forall f \in E \setminus S(X,Y): Z_f = (X \circ Y)_f.
\end{align*} 
\end{itemize}
The elements of $\mathcal{L}$ are called \emph{covectors}.
\end{Definition}

In order to see that the zero sets of a COM form a bouquet of geometric lattices, we need to define oriented matroids first. The most common definition states that $(E,\mathcal{L})$ forms an \emph{oriented matroid (OM)}, if $\mathcal{L}$ fulfills

\begin{itemize}
    \item[(S)]  Symmetry 
    \begin{align*}
\forall X\in \mathcal{L}: -X \in \mathcal{L}.
\end{align*}
\end{itemize}
and
\begin{itemize}
\item[(C)] Composition
\begin{align*}
\forall X,Y \in \mathcal{L}: X \circ (-Y) \in \mathcal{L}.
\end{align*} 
\end{itemize}

and (SE). It is easy to check that (FS) implies (C). So in order to be an OM, a COM just has to fullfil (S) additionally. Note that if a COM contain the all zero vector $\mathbf{0}$, (S) is fulfilled because of (FS). So we can define OMs in terms of COMs

\begin{Definition}[Oriented Matroid (OM)] Let $E$ be a finite set and $\mathcal{L} \subseteq \{0,+,-\}^{|E|}$. The pair $\mathcal{M}=(E,\mathcal{L})$ is called an OM, if it is a COM that contains the all zero vector $\mathbf{0}$.
\end{Definition}

Let $\mathcal{M} = (E,\mathcal{L})$ be a COM. For a covector $X\in \mathcal{L}$, the set $F(X)=\{X\circ Y\mid Y\in\mathcal{L}\}$ is called the \emph{face} of $X$. Note that $(E\backslash\underline{X}, F(X)\backslash \underline{X})$ is an oriented matroid. This can be checked by looking at the definition of OMs given above. With this observation we can use the following well known result.

\begin{Theorem}[\cite{bjorner1999oriented}]\label{thm:OMzero}
  Let $\mathcal{L}=\{+, -, 0\}^{|E|}$ be the set of covectors of an oriented matroid. Then the set $\{z(X):\,X \in \mathcal{L}\}$ is the collection of flats of the underlying matroid.
\end{Theorem}

Using that result it is easy to proof the following theorem.

\begin{Theorem}
    Let $(E,\mathcal{L})$ be a COM. The poset $Z$ on $\{z(X):\,X \in \mathcal{L}\}$ induced by inclusion is a bouquet of geometric lattices.
\end{Theorem}

\begin{proof}
Accourding to \ref{def:bouquetlattice} we have to check two things, first that $Z$ is a meet semilattice and second that every interval is a geometric lattice. To see that $Z$ is a meet semilattice let $X,Y \in \mathcal{L}$. For the two elements $z(X)$ and $z(Y)$ of $Z$, $z(X \circ Y)$ is a meet. So every two elements have a meet, which shows that $Z$ is a meet semilattice. Now let $\{z(X),z(Y)\}$ be an interval in $Z$. We know that $(E\backslash\underline{Y}, F(Y)\backslash \underline{Y})$ is an OM. Using \ref{thm:OMzero} we see that everything that lies below $z(Y)$ is a geometric lattice. Since every interval of a geometric lattice is again geometric, $\{z(X),z(Y)\}$ is geometric as well.
\end{proof}

Now we can apply \ref{main} and get the following corollary.

\begin{Corollary}
Let $(E,\mathcal{L})$ be a COM, $Z$ the poset on $\{z(X):\,X \in \mathcal{L}\}$ induced by inclusion and $\mathfrak{C}$ its chain matrix. Then 
\begin{align*}
    det(\mathfrak{C}) = \pm \prod_{Y \in \mathcal{L}} w(Y)^{\rho_P(x)}.
\end{align*}
where $w(x) = \sum_{e \in z(Y)}w_e$.
\end{Corollary}

\subsection{Bouquets of matroids}\label{sec:bouquetsMatr}
Like matroids, bouquets of matroids can be defined in various ways. To stay consistent we choose the definition by independent sets.

\begin{Definition}[\cite{laurent1989bouquets}]
Let $E_1,\dots,E_n$ be a clutter of subsets of $E$ and $\mathcal{I}$ a collection of subsets of $E$. $\mathcal{I}$ is called a bouquet of matroids on $E$ with roofs $E_1,\dots,E_n$ if
\begin{enumerate}
    \item $\mathcal{I}_i = \mathcal{I}\cap 2^{X_i}$ is the family of independent sets of a matroid on $E_i$ for all $i = 1,\dots,n$. 
    \item $\mathcal{I} = \bigcup^m_{i=1} \mathcal{I}_i$
    \item If $I \in \mathcal{I}_i \cap \mathcal{I}_j$ and $e \in E_i - E_j$, then $I\cup e \in \mathcal{I}$.
\end{enumerate}
\end{Definition}
The definition of flats and flags of a bouquet of matroids is analog to the definition of flats and flags of a matroids. The following gives the connection to bouquets of geometric lattices:

\begin{Theorem}[\cite{laurent1989bouquets}]
The poset of flats of a bouquet of matroids is a bouquet of geometric lattices.
\end{Theorem}

We say that the chain matrix of a bouquet of matroids is called the flag matrix, like in the matroid case. From \ref{main} is follows immediately
\begin{Corollary}
Let $\mathcal{M}$ be a bouquet of matroids and $\mathfrak{F}$ its flag matrix. Then 
\begin{align*}
det(\mathfrak{F}) = \pm \prod\limits_{K \text{ flat of } \mathcal{M}} w(K)^{\rho_P(x)}
\end{align*}
where $w(K) = \sum_{e \in K}w_e$.
\end{Corollary}

\bibliographystyle{plain} 
\bibliography{BVbib}

\begin{thebibliography}{10}

\bibitem{bandelt2018coms}
Hans-J{\"u}rgen Bandelt, Victor Chepoi, and Kolja Knauer.
\newblock Coms: complexes of oriented matroids.
\newblock {\em Journal of Combinatorial Theory, Series A}, 156:195--237, 2018.

\bibitem{bjorner1999oriented}
Anders Bj{\"o}rner.
\newblock {\em Oriented matroids}.
\newblock Number~46. Cambridge University Press, 1999.

\bibitem{brylawski1997determinant}
T~Brylawski and A~Varchenko.
\newblock The determinant formula for a matroid bilinear form.
\newblock {\em advances in mathematics}, 129(1):1--24, 1997.

\bibitem{brylawski2000mobius}
Tom Brylawski.
\newblock A m{\"o}bius identity arising from modularity in a matroid bilinear
  form.
\newblock {\em Journal of Combinatorial Theory, Series A}, 91(1-2):622--639,
  2000.

\bibitem{crapo1967higher}
Henry~H Crapo.
\newblock A higher invariant for matroids.
\newblock {\em Journal of Combinatorial Theory}, 2(4):406--417, 1967.

\bibitem{eberhardt2024intersection}
Jens~Niklas Eberhardt and Carl Mautner.
\newblock An intersection matrix for affine hyperplane arrangements.
\newblock {\em arXiv preprint arXiv:2407.06008}, 2024.

\bibitem{hochstattler2022signed}
Winfried Hochst{\"a}ttler, Sophia Keip, and Kolja Knauer.
\newblock The signed varchenko determinant for complexes of oriented matroids.
\newblock {\em arXiv preprint arXiv:2211.13986}, 2022.

\bibitem{hochstattler2019varchenko}
Winfried Hochst{\"a}ttler and Volkmar Welker.
\newblock The varchenko determinant for oriented matroids.
\newblock {\em Mathematische Zeitschrift}, 293:1415--1430, 2019.

\bibitem{laurent1989bouquets}
Monique Laurent and Michel Deza.
\newblock Bouquets of geometric lattices: some algebraic and topological
  aspects.
\newblock In {\em Annals of Discrete Mathematics}, volume~43, pages 279--313.
  Elsevier, 1989.

\bibitem{schechtman1991arrangements}
Vadim~V Schechtman and Alexander~N Varchenko.
\newblock Arrangements of hyperplanes and lie algebra homology.
\newblock {\em Inventiones mathematicae}, 106(1):139--194, 1991.

\bibitem{varchenko1995multidimensional}
Alexander Varchenko.
\newblock {\em Multidimensional hypergeometric functions the representation
  theory of Lie Algebras and quantum groups}, volume~21.
\newblock World Scientific, 1995.

\bibitem{wachs2006poset}
Michelle~L Wachs.
\newblock Poset topology: tools and applications.
\newblock {\em arXiv preprint math/0602226}, 2006.

\bibitem{wachs1985geometric}
Michelle~L Wachs and James~W Walker.
\newblock On geometric semilattices.
\newblock {\em Order}, 2(4):367--385, 1985.

\end{thebibliography}

\end{document}